\documentclass[11pt, reqno]{amsart}

\usepackage{amssymb}  % AMS fonts

\theoremstyle{plain}
\newtheorem{theorem}{Theorem}[section]
\newtheorem{lemma}[theorem]{Lemma}
\newtheorem{corollary}[theorem]{Corollary}
\newtheorem{proposition}[theorem]{Proposition}
\newtheorem{conjecture}[theorem]{Conjecture}
\newtheorem{problem}[theorem]{Problem}

\theoremstyle{definition}

\theoremstyle{remark}
\newtheorem{remark}[theorem]{Remark}
\newtheorem{example}[theorem]{Example}

\newcommand{\seclabel}[1]{\label{sec:#1}} % section
\newcommand{\subseclabel}[1]{\label{subsec:#1}} % subsection
\newcommand{\thmlabel}[1]{\label{thm:#1}} % theorem
\newcommand{\lemlabel}[1]{\label{lem:#1}} % lemma
\newcommand{\corolabel}[1]{\label{coro:#1}} % corollary
\newcommand{\proplabel}[1]{\label{prop:#1}} % proposition
\newcommand{\remlabel}[1]{\label{rem:#1}} % remark
\newcommand{\problabel}[1]{\label{prob:#1}} % problem
\newcommand{\tablelabel}[1]{\label{table:#1}} % table
\newcommand{\eqlabel}[1]{\label{eq:#1}} % equation

\newcommand{\secref}[1]{\ref{sec:#1}} % section
\newcommand{\thmref}[1]{\ref{thm:#1}} % theorem
\newcommand{\lemref}[1]{\ref{lem:#1}} % lemma
\newcommand{\cororef}[1]{\ref{coro:#1}} % corollary
\newcommand{\propref}[1]{\ref{prop:#1}} % proposition
\newcommand{\probref}[1]{\ref{prob:#1}} % problem
\newcommand{\tableref}[1]{\ref{table:#1}} % table
\renewcommand{\eqref}[1]{\ref{eq:#1}} % no parentheses
\newcommand{\peqref}[1]{(\eqref{#1})} % parenthesized equation ref

\newcommand{\aut}[1]{\mathrm{Aut}\,#1}
\newcommand{\ext}[2]{\ltimes_{#1}^{#2}}
\newcommand{\im}[1]{\mathrm{Im}\,#1}
\renewcommand{\ker}[1]{\mathrm{Ker}\,#1}
\newcommand{\ov}[1]{\overline{#1}}

\title[C-loops: Extensions and Constructions]
{C-loops: Extensions and Constructions}

\author[M.~K.~Kinyon]{Michael~K.~Kinyon}
\address{Department of Mathematical Sciences \\
Indiana University South Bend \\
South Bend, IN 46634 USA}
\email{mkinyon@iusb.edu}
\urladdr{http://mypage.iusb.edu/\symbol{126}mkinyon}

\author[J.~D.~Phillips]{J.~D.~Phillips}
\address{Department of Mathematics \& Computer Science \\
Wabash College \\
Crawfordsville, IN 47933 U.S.A.}
\email{phillipj@wabash.edu}
\urladdr{http://www.wabash.edu/depart/math/faculty.html{\#}Phillips}

\author[P. Vojt\v{e}chovsk\'{y}]{Petr Vojt\v{e}chovsk\'{y}}
\address{Department of Mathematics \\
University of Denver \\
2360 S Gaylord St \\
Denver, CO 80208 U.S.A.}
\email{petr@math.du.edu}
\urladdr{http://www.math.du.edu/\symbol{126}petr}

\date{\today}


\begin{document}

\begin{abstract}
C-loops are loops satisfying the identity $x(y\cdot yz) = (xy\cdot y)z$.
We develop the theory of extensions of C-loops, and characterize all
nuclear extensions provided the nucleus is an abelian group. C-loops with central
squares have very transparent extensions; they can be built from small blocks
arising from the underlying Steiner triple system. Using these extensions,
we decide for which abelian groups $K$ and Steiner loops $Q$ there is a
nonflexible C-loop $C$ with center $K$ such that $C/K$ is isomorphic to $Q$.
We discuss possible orders of associators in C-loops. Finally, we show that
the loops of signed basis elements in the standard real Cayley-Dickson algebras
are C-loops.
\end{abstract}

\maketitle

\section{Introduction}
\seclabel{intro}

\noindent This is the second paper in a series devoted to \emph{C-loops},
which are loops satisfying the identity $x(y\cdot yz) = (xy\cdot y)z$.
C-loops were introduced by Fenyves \cite{Fe2}. Moufang
C-loops are exactly Fenyves' \emph{extra loops} \cite{Fe1, extra}.
The first detailed study of general C-loops was in
\cite{CLoopsI}. In particular, it was observed that
C-loops can be characterized as inverse property loops
with all squares in the nucleus.
We also examined basic loop-theoretical properties of C-loops,
established some connections between C-loops, Steiner loops
and commutative Moufang loops, proved a structural result for torsion commutative C-loops, and obtained the first few smallest nonextra C-loops.

In the present paper, we are concerned with constructions of C-loops.
We assume some familiarity with \cite{CLoopsI}, but results from
there are quoted as needed. For general background in loop theory,
the standard references are \cite{bel,bruck,cps,pflug}.

In \S\secref{ext}, we develop the theory of extensions for C-loops.
Given an abelian group $K$ and a C-loop $Q$ we show that a C-loop $C$ is an
extension of $K$ by $Q$ if and only if there is a C-factor set (C-cocycle)
satisfying a certain condition; see Theorem \thmref{ext}.

In \S\secref{some} we consider central extensions, that is, when
$K$ is contained in the center $Z(C)$ of the C-loop extension $C$.
For C-loops with central squares, the
C-factor set condition for central extensions is greatly simplified.
Moreover, the corresponding extensions
can be constructed from small blocks arising from the combinatorial
properties of the underlying Steiner triple system. These
``block extensions'' are described in \S\secref{BB}.

If $L$ is a C-loop with nucleus $N=N(L)$, then $L/N$ is a Steiner loop.
Furthermore, when $Z=Z(L)$ is the center of $L$ then $L/Z$ is Steiner if and
only if every square in $L$ is central. One of the main themes of the paper is
the following question:

\noindent\emph{For which abelian groups $K$ and Steiner loops $Q$ is
there a C-loop $C$ with nucleus/center isomorphic to $K$ such that
$C/K$ is isomorphic to $Q$?}

\noindent A loop
is \emph{flexible} if it satisfies the identity $x\cdot yx = xy\cdot x$.
In particular, commutative loops are flexible.  
Flexible C-loops are diassociative, \emph{i.e.}, any two elements 
generate an associative subloop \cite[Lemma 4.4]{CLoopsI}. It turns
out to be easier to answer our question in the \emph{non}flexible case.
Besides, focusing on nonflexible C-loops ensures that the extension
is neither Moufang nor Steiner.

When the abelian group $K$ is the center of the nonflexible C-loop $C$,
we answer our question completely (cf. Theorem \thmref{NN}). When $K$
is the nucleus of $C$, the only unresolved case is where $K$ is an
elementary abelian $2$-group, and even then we have partial results
(cf. \S\secref{NN}).

We conclude \S\secref{NN} with a construction that yields nonflexible
C-loops that possess an associator of given order $n$, where $n > 2$.

Next, consider the standard Cayley-Dickson algebras over the real numbers.
It is possible to choose their bases in such a way that the signed basis
elements form a loop under multiplication regardless of the dimension.
This is shown in \S\secref{cd}. In fact, the loops are flexible C-loops,
and if the dimension of the algebra is at least $32$, then the loops are
neither Moufang nor commutative, and have nucleus (and hence center)
of order $2$; see \S\secref{sed}.

Conjectures and open problems are presented throughout the paper.

\section{Nuclear extensions}
\seclabel{ext}

For the general extension theory of loops, we refer the reader to
\cite[Chapter III]{cps}. Here we limit ourselves to C-loops.
Our first goal is to characterize a broad class of extensions of
C-loops. Our notation for extensions is set up to resemble that of
\cite{Rotman}.

Let $K$, $Q$ be C-loops. Then a loop $C$ is said to be an \emph{extension of
$K$ by $Q$} if $K$ is a normal subloop of $C$ such that $C/K$ is isomorphic
to $Q$.

In the above situation, let $\pi:C\to C/K=Q$ be the natural projection, and
let $\ell:Q\to C$ be a \emph{section} of $\pi$, i.e, $\pi\ell x = x$
for every $x\in Q$. Throughout the paper, we assume that $\ell 1=1$
whenever $\ell$ is a section. Then, for $x,y\in Q$, we have
$K(\ell x\cdot \ell y) = K(\ell(xy))$, and there is therefore a
unique element $f(x,y)\in K$ such that
\[
    \ell x\cdot \ell y = f(x,y)\cdot \ell(xy).
\]
The resulting map $f:Q\times Q\to K$ is said to be \emph{associated with}
$\ell$.

Let $\theta:Q\to\aut{K}; x\mapsto \theta_x$ be a homomorphism. The
pair $(\theta,f)$ is said to be a \emph{C-factor set}
(or \emph{C-cocycle}) if
\begin{equation}
\eqlabel{cfs1}
    f(x,1)=1=f(1,x),
\end{equation}
and
\begin{equation}
\eqlabel{cfs2}
    \theta_{xy}f(y,z)\cdot \theta_x f(y,yz)\cdot f(x,y\cdot yz)
    = f(x,y)\cdot f(xy,y) \cdot f(xy\cdot y,z)
\end{equation}
holds for every $x$, $y$, $z\in Q$.

\begin{remark}
\remlabel{ambiguous}
Equation $\peqref{cfs2}$ is ambiguous unless $K$ is associative.
In all situations discussed below, $K$ is an abelian group. A less general
notion of a C-factor set $($with trivial $\theta)$ was introduced in
\cite{CLoopsI}.
\end{remark}

Given a C-factor set $(\theta,f)$, we define a binary operation
$*$ on $D = K\times Q$ by
\begin{equation}
 \eqlabel{KQ}
    (a,x)*(b,y) = (a\cdot\theta_xb\cdot f(x,y),\, xy),
\end{equation}
for $a, b\in K$ and $x, y\in Q$. We denote the resulting
quasigroup by $K\ext{\theta}{f}Q$.

For a loop $C$ and $x\in C$, we denote the left and right
multiplication maps by $L_x : C\to C; y\mapsto xy$ and
$R_x:C\to C; y\mapsto yx$, respectively. Set $T_x=R_x^{-1}L_x$;
these are usually called \emph{middle inner mappings}.
In group theory, $T_x$ is conjugation by $x$, and is an
automorphism. This is not necessarily so when $C$ is not
associative. Roughly speaking, the theory of extensions of
loops can imitate the theory of extensions of groups so long
as the middle inner mappings $T_x$ behave as automorphisms
on the normal subloops in question, and satisfy $T_xT_y=T_{xy}$.

Since every C-loop is an inverse property loop
(\cite[Corollary 2.4]{CLoopsI}), we have $T_x=R_{x^{-1}}L_x$.

\begin{proposition}
\proplabel{D}
Let $K$ be an abelian group, $Q$ a C-loop, $f:Q\times Q\to K$ a map, and
$\theta: Q\to \aut{K}$ a homomorphism. Then $D=K\ext{\theta}{f} Q$ is a
C-loop with neutral element $(1,1)$ if and only if $(\theta,f)$ is a C-factor
set. Moreover, when $D$ is a C-loop, we have:
\begin{enumerate}
\item[(i)] $K\cong (K,1)\le D$, $Q\cong (1,Q)\le D$,

\item[(ii)] $K\le N(D)$,

\item[(iii)] for every $a\in K$, $x\in Q$,
$(a,x)^{-1}=(\theta_x(f(x,x^{-1})a),\,x^{-1})$,

\item[(iv)] for every $a\in K$, $x\in Q$, $\theta_x=T_{(a,x)}|_K$,

\item[(v)] $K\unlhd D$, $D/K\cong Q$.
\end{enumerate}
\end{proposition}
\begin{proof}
Let $u=(a,x)$, $v=(b,y)$, $w=(c,z)$, for $a$, $b$, $c\in K$, $x$, $y$, $z\in
Q$. Then
\begin{align*}
    u*(v*(v*w))&=u*(v*(b\theta_ycf(y,z),\,yz))\\
    &=u*(b\theta_y(b\theta_ycf(y,z))f(y,yz),\,y\cdot yz)\\
    &=(a\theta_x(b\theta_y(b\theta_ycf(y,z))f(y,yz))f(x,y\cdot yz),\,
        x(y\cdot yz))\\
    &=(a\theta_xb\theta_{xy}b\theta_{xy\cdot y}c\theta_{xy}
        f(y,z)\theta_{x}f(y,yz)f(x,y\cdot yz),\, x(y\cdot yz)).
\end{align*}
On the other hand
\begin{align*}
    ((u*v)*v)*w&=((a\theta_xbf(x,y),\,xy)*v)*w\\
    &=(a\theta_xbf(x,y)\theta_{xy}bf(xy,y),\,xy\cdot y)*w\\
    &=(a\theta_xbf(x,y)\theta_{xy}bf(xy,y)\theta_{xy\cdot y}c
        f(xy\cdot y,z),\,(xy\cdot y)z).
\end{align*}
Since $Q$ is a C-loop and $K$ is an abelian group, we see that
$u*(v*(v*w))=((u*v)*v)*w$ if and only if \peqref{cfs2} holds.

We have $(a,x)*(1,1)=(af(x,1),x)$ and $(1,1)*(a,x)=(af(1,x),x)$. Thus $(1,1)$
is the neutral element of $D$ if and only if \peqref{cfs1} holds.

For the rest of the proof, assume that $D$ is a C-loop. Part (i) is
straightforward. For $a,b,c\in K$, $y,z\in Q$, we have
$(a,1)*((b,y)*(c,z)) = (a,1)*(b\theta_ycf(y,z),\,yz) =
(ab\theta_ycf(y,z),\,yz) = (ab,y)*(c,z) = ((a,1)*(b,y))*(c,z)$. Therefore $K$
is contained in the left nucleus of $D$. By \cite[Corollary 2.5]{CLoopsI},
the three nuclei of $D$ coincide, and (ii) follows.

As inverse property loops, C-loops satisfy the identity
$(xy)^{-1}=y^{-1}x^{-1}$. Let $(a,x)\in D$. Then
\begin{align*}
    &(a,\,x)(\theta_x(f(x,x^{-1})a),\,x^{-1}) =
    (a\theta_x\theta_x(f(x,x^{-1})a)f(x,x^{-1}),\,1)\\
    &=(a\theta_x\theta_x^{-1}(a^{-1}f(x,x^{-1})^{-1})f(x,x^{-1}),\,1)
    =(1,\,1),
\end{align*}
proving (iii).

Upon identifying $K$ with $(K,1)$ and $Q$ with $(1,Q)$, it makes sense to
write $\theta_x(k,1) = \theta_xk = (\theta_xk,1)$ for every $k\in K$, $x\in
Q$. Given $(a,x)\in D$, a short calculation yields $T_{(a,x)}(k,1) =
(a,x)(k,1)\cdot (a,x)^{-1} = (\theta_xk,1)$, and we are done with (iv).

In order to check that $K$ is normal in $D$, we must show that $K$ is 
invariant under the standard generators for the inner mapping group: $L(x,y)=L_{yx}^{-1}L_yL_x$, $R(x,y)=R_{xy}^{-1}R_yR_x$, $T_x$, for
$x, y\in D$. Since $K$ is nuclear, $K$ is trivially invariant under
each $L(x,y)$ and $R(x,y)$. By (iv), $K$ is invariant under each $T_x$,
and so $K\unlhd D$ follows. It is then easy to show that $K/D$ is isomorphic
to $Q$.
\end{proof}

The following result of Leong \cite[Theorem 3]{Leong} is important here:

\begin{lemma}\lemlabel{Leong}
Let $Q$ be a loop with a normal subloop $K \leq N(Q)$. For each $x \in Q$,
define $\theta_x = T_x |_K$. Then
\begin{enumerate}
\item For each $x\in Q$, $\theta_x \in Aut(K)$,

\item The mapping $\theta : Q\to Aut(K)$ is a homomorphism.
\end{enumerate}
\end{lemma}
\begin{proof}
First fix $a, b\in K$ and $x\in Q$. Since $K$ is normal in $Q$, we have
$\theta_x(K)\le K\le N(Q)$. In particular,
$T_x(ab)\cdot x = x\cdot ab = xa\cdot b = (T_x(a)\cdot x)b
= T_x(a)\cdot xb = T_x(a) (T_x(b)\cdot x) = T_x(a) T_x(b) \cdot x$.
Cancelling $x$ on the right then shows that
$\theta_x$ is an automorphism of $K$.

Now fix $a\in K$ and $x,y\in Q$. Let $z = T_{xy}a$. By the
first part, $z\in K$, and so $(zx)y = z(xy) = (xy)a = x(ya) =
x (T_y(a)\cdot y) = x T_y(a)\cdot y$. Upon cancelling $y$
on the right, we get $zx = x T_y(a)$, i.e., $z = T_x T_y (a)$.
Hence $\theta_{xy}=\theta_x\theta_y$ as claimed.
\end{proof}

\begin{theorem}
\thmlabel{ext}
Let $K$ be an abelian group and $Q$ a C-loop. The following conditions are
equivalent:
\begin{enumerate}
\item[(i)] $C$ is a C-loop extension of $K$ by $Q$, and $K\le N(C)$,

\item[(ii)] there is a homomorphism $\theta:Q\to\aut{K}$ and a map $f:Q\times
Q\to K$ such that $(\theta,f)$ is a C-factor set and $C=K\ext{\theta}{f}Q$.
\end{enumerate}
\end{theorem}

\begin{proof}
If (ii) holds, (i) follows by Proposition \propref{D}.

Conversely, assume that (i) holds. Let $\pi:C\to C/K=Q$ be the natural
projection, and let $\ell:Q\to C$ be a section. Let $f:Q\times Q\to K$ be
associated with $\ell$. By Lemma \lemref{Leong}, the map $\theta:Q\to
\aut{K}$, $x\mapsto \theta_x=T_{\ell(x)}|_K$ is well-defined. Given
$x, y\in Q$, let $n=f(x,y)\in N(C)$. Then, by Lemma \lemref{Leong},
$\theta_{x}\theta_{y}=T_{n\ell(xy)} = T_n\theta_{xy} = \theta_{xy}$, showing
that $\theta$ is a homomorphism. We proceed to show that $(\theta,f)$ is a
C-factor set.

Let $x, y, z\in Q$. On the one hand,
\begin{align*}
    &(\ell(x)\ell(y)\cdot \ell(y))\ell(z)
        =(f(x,y)\ell(xy)\cdot \ell(y))\ell(z)
        = f(x,y)(\ell(xy)\ell(y)\cdot\ell(z))\\
    &=f(x,y)f(xy,y)\ell(xy\cdot y)\ell(z)
        =f(x,y)f(xy,y)f(xy\cdot y,z)\ell((xy\cdot y)z),
\end{align*}
as $K\le N(C)$. On the other hand, since $C$ is an inverse property loop,
$K\le N(C)$ and $T_{\ell(x)}|_K=\theta_x\in \aut{K}$, we have
\[
    \ell(x)\cdot uv = \ell(x)u\cdot v
    = (T_{\ell(x)}u\cdot \ell(x)) v
    = (\theta_x u\cdot \ell(x))v = \theta_xu\cdot
    \ell(x)v
\]
 for every $u\in K$, $v\in C$, and therefore
\begin{align*}
    &\ell(x)(\ell(y)\cdot\ell(y)\ell(z))
        = \ell(x)(\ell(y)\cdot f(y,z)\ell(yz))
        = \ell(x)(\theta_yf(y,z)\cdot \ell(y)\ell(yz))\\
    &=\ell(x)(\theta_yf(y,z)\cdot f(y,yz)\ell(y\cdot yz))
        = \ell(x)\theta_yf(y,z)\cdot f(y,yz)\ell(y\cdot yz)\\
    &=\theta_x\theta_yf(y,z)\cdot \ell(x)(f(y,yz)\ell(y\cdot yz))
       = \theta_{xy}f(y,z)\cdot \theta_xf(y,yz)\cdot \ell(x)\ell(y\cdot yz)\\
    &= \theta_{xy}f(y,z)\cdot \theta_xf(y,yz)\cdot f(x,y\cdot yz)\cdot
        \ell(x,y\cdot yz).
\end{align*}
As both $C$ and $Q$ are C-loops, we deduce that $(\theta,f)$ satisfies
\peqref{cfs2}. Since \peqref{cfs1} holds by definition of $f$,
$(\theta,f)$ is a C-factor set.

It remains to show that there is an isomorphism $\psi$ from $C$ to
$(D,*)=K\ext{\theta}{f}Q$. Given $u\in C$, there are uniquely determined
elements $a\in K$, $x\in Q$ such that $u=ax$. Accordingly, we set
$\psi(u)=(a,x)$. Then $\psi$ is a bijection, and it suffices to prove that it
is a homomorphism. Consider $v=by$ where $b\in K$, $y\in Q$. Since $a, b\in
N(C)$, we have $u\cdot v = a(xb\cdot y)$, which can be rewritten as
$a(\theta_xb\cdot xy) = a\theta_xb\cdot  xy$, using the same trick as above.
Then $\psi(u\cdot v) = (a\theta_xb,\,xy) = (a,x)*(b,y) = \psi(u)*\psi(v)$.
\end{proof}

Since we would like to understand the extensions of C-loops up to
isomorphism, we ask:

\begin{problem}
Let $K$ be an abelian group, $Q$ a C-loop, and $(\theta, f)$, $(\theta',f')$
two C-factor sets. Under which conditions on $(\theta,f)$, $(\theta',f')$ are
the C-loops $K\ext{\theta}{f}Q$, $K\ext{\theta'}{f'}Q$ isomorphic?
\end{problem}

\subsection{C-loops with nonabelian nucleus}
\subseclabel{noncomm-nuc}

\noindent Theorem \thmref{ext} does not capture all nuclear extensions for
C-loops since there are C-loops with nonabelian nucleus. For instance,
let $D$ be the direct product of the symmetric group $S_3$ and the smallest
nonassociative Steiner loop $Q$ of order $10$. Then $D$ is a C-loop of
order $60$, and the nucleus of $D$ contains (and in fact, coincides with)
$S_3$, a nonabelian group.

Recall the following result:

\begin{theorem}
\thmlabel{Folk}
Let $L$ be a nonassociative loop with normal nucleus $N = N(L)$
such that $L/N$ is a group. Then $N$ has nontrivial center.
\end{theorem}

\begin{proof}
This can be found in, for instance, \cite[\S4]{KKP}.
\end{proof}

\begin{corollary}
\corolabel{ten}
Let $C$ be a nonassociative C-loop with nucleus $N = N(C)$,
and assume that $N$ has trivial center. Then $C/N$ is a
nonassociative Steiner loop. In particular, $|C|\ge 10|N|$.
\end{corollary}

\begin{proof}
The quotient $C/N$ is a Steiner loop by \cite[Proposition 5.8]{CLoopsI}.
Since $C$ is nonassociative, Theorem \thmref{Folk} implies that
$C/N$ cannot be a group. The smallest nonassociative Steiner loop
has order $10$.
\end{proof}

This shows that the above direct product is the smallest nonassociative
C-loop with nucleus isomorphic to $S_3$. However, it is conceivable that
there is a smaller C-loop with nonabelian nucleus, provided the nucleus
does not have trivial center. We therefore ask:

\begin{problem} What is the smallest integer $n$ for which there exists a
C-loop of order $n$ with nonabelian nucleus?
\end{problem}

Exhaustive computer searches suggest that there is no C-loop
of order $32$ with a nucleus containing a nonabelian group of
order $8$.

\section{Central extensions}
\seclabel{some}

The C-loop extensions described in the previous section can, as the
section's title suggests, be considered to be nuclear extensions.
Not surprisingly, central C-loop extensions are characterized
by the triviality of the homomorphism $\theta$.

\begin{lemma}
\lemlabel{Z}
Let $K$ be an abelian group, $Q$ a C-loop, $(\theta,f)$ a C-factor set and
$C=K\ext{\theta}{f}Q$. Then $K\le Z(C)$ if and only if $\theta=\mathrm{id}$.
\end{lemma}

\begin{proof}
If $\theta=\mathrm{id}$, then for every $a, b\in K$ and $x\in Q$,
we have $(a,1)(b,x)=(ab,x)=(ba,x)=(b,x)(a,1)$, so that $K\le Z(C)$.

Conversely, if $K\le Z(C)$, then $(a,1)(b,x)=(ab,x)$ is equal to
$(b,x)(a,1)=(b\theta_x(a),x)$ for every $a, b\in K$ and $x\in Q$. This
means that $x=\theta_x(a)$ for every $a\in K$, $x\in Q$.
\end{proof}

Let $C$ be a C-loop with nucleus $N = N(C)$. Then $N$ is
normal in $C$ and $C/N$ is a Steiner loop. It is not true, however,
that $C/Z(C)$ is necessarily a Steiner loop. Recalling that Steiner
loops are precisely inverse property loops of exponent two
(\cite[Lemma 2.2]{CLoopsI}), we immediately obtain:

\begin{proposition}
\proplabel{steiner-squares}
Let $L$ be a C-loop. Then $L/Z(L)$ is a Steiner loop if and only if all
squares of $L$ are central.
\end{proposition}

Next we consider the situation where the Steiner quotient $C/N$ is simple.
Note that we do not assume that $C$ is a C-loop in the following
lemma:

\begin{lemma}
\lemlabel{sim}
Let $C$ be a loop with abelian nucleus $N(C)$. If $C/N(C)$ is
simple and nonassociative, then $N(C)=Z(C)$.
\end{lemma}

\begin{proof}
Let $N=N(C)$. Consider again the homomorphism $\theta:C\to \aut{N};
x\mapsto \theta_x= T_x|_N$. As in the proof of Theorem \thmref{ext},
$\theta_n\theta_x=\theta_x$ for every $n\in N$. Then $\theta$ induces a
homomorphism $\ov{\theta}:C/N\to \aut{N}$,
$\ov{\theta}_x=\theta_x$. Assume that $\ov{\theta}$ is injective.
Then $\im{\ov{\theta}}\simeq C/N$ is nonassociative; a contradiction
with $\im{\ov{\theta}}\le \aut{N}$. Thus $\ov{\theta}$ is not
injective. Since $C/N$ is simple, we conclude that
$\ker{\ov{\theta}}=C/N$, i.e., $\theta_x=\mathrm{id}_N$ for every $x\in
C$. Then $N(C)\le Z(C)$ follows.
\end{proof}

\begin{proposition}
\proplabel{sim}
Let $C$ be a C-loop with abelian nucleus $N(C)$. If $C/N(C)$ is
simple, then $N(C)=Z(C)$.
\end{proposition}

\begin{proof}
By Lemma \lemref{sim}, we may assume $Q=C/N(C)$ is a simple,
associative Steiner loop. Since Steiner loops are
commutative of exponent $2$, $Q$ is either trivial or a cyclic group
of order $2$. In the former case, $C=N(C)$ is an abelian group and $Z(C)=N(C)$
follows. The latter case never occurs, since no nonassociative C-loop has
nucleus of exponent $2$, by \cite[Lemma 2.9]{CLoopsI}.
\end{proof}

By a result of Quackenbush \cite{Qu}, given a Steiner triple system,
either its associated Steiner quasigroup or its associated Steiner loop
is simple. The hypotheses of Proposition \propref{sim}
are therefore often fulfilled when $C$ is a C-loop.

For the rest of the paper, we will consider only the case
that the quotient $C/Z(C)$ is a Steiner loop.

The definition of a C-factor set can be greatly simplified when
$\theta=\mathrm{id}$ and when $Q$ is a Steiner loop. 
Except as otherwise noted, we will also write the abelian
group $K$ additively.

\begin{proposition} Let $Q$ be a Steiner loop, $K$ an abelian group and
$f:Q\times Q\to K$ a map satisfying \peqref{cfs1}. Then the following
conditions are equivalent:
\begin{enumerate}
\item[(i)] $(\mathrm{id},f)$ is a C-factor set,

\item[(ii)] for every $x,y,z\in Q$, $f$ satisfies
\begin{equation} \eqlabel{CFS2}
    f(y,z)+f(y,yz)=f(x,y)+f(xy,y),
\end{equation}

\item[(iii)] for every $x,y,z\in Q$, $f$ satisfies the two conditions
\begin{align}
    f(xy,y)&=f(y,y)-f(x,y), \eqlabel{MK1}\\
    f(y,yz)&=f(y,y)-f(y,z). \eqlabel{MK2}
\end{align}
\end{enumerate}
\end{proposition}

\begin{proof}
Since $Q$ satisfies $y\cdot yz=z$ and $xy\cdot y=x$, \peqref{CFS2}
is equivalent to \peqref{cfs2}. Assume that \peqref{CFS2} holds. Then
\peqref{MK1} is obtained from \peqref{CFS2} with $z=1$, and
\peqref{MK2} is obtained from \peqref{CFS2} with $x=1$. Conversely,
assume that \peqref{MK1}, \peqref{MK2} hold. Then subtracting
\peqref{MK1} from \peqref{MK2} gives \peqref{CFS2}.
\end{proof}

The advantage of equations \peqref{MK1}, \peqref{MK2} over
\peqref{CFS2} is that they deal with only $2$ elements at a time, and are
therefore easier to verify. By imposing another condition on $f$, the
equations become even simpler.

\begin{lemma}\lemlabel{Aux}
Let $Q$ be a Steiner loop, $K$ an abelian group, and $f:Q\times Q\to
K$ a map satisfying \peqref{cfs1} and $f(y,y)=0$ for every $y\in L$. Then
$K\ext{\mathrm{id}}{f}Q$ is a C-loop if and only if both
\begin{align}
    f(xy,y)&=-f(x,y), \eqlabel{Q1}\\
    f(y,yz)&=-f(y,z) \eqlabel{Q2}
\end{align}
hold for every $x,y,z\in L$.
\end{lemma}

\section{The building blocks of central extensions}
\seclabel{BB}

\noindent We now build all central extensions with Steiner quotient $Q$. We
construct all of them from small pieces arising from the underlying Steiner
triple system.

In a sense, all Steiner loops are locally the same. Consider a Steiner loop
$Q$ and $x,y\in Q$ with $1\ne x\ne y\ne 1$. Since $\{x,y\}$ is one of the
edges of the underlying Steiner triple system, the subloop $\langle
x,y\rangle$ generated by $\{x,y\}$ corresponds to one triangle (block) of the
Steiner triple system, and is isomorphic to the Klein group.

Since the defining equations \peqref{MK1}, \peqref{MK2} of a central
C-factor set for $Q$ Steiner deal only with two elements at a time, the map
$f$ can be build by small pieces:

\begin{proposition} Let $Q$ be a Steiner loop, $K$ an abelian group, and
$f:Q\times Q\to K$ a map. Then $(\mathrm{id},f)$ is a C-factor set if and
only if for every $u$, $v\in Q$ with $1\ne u\ne v\ne 1$ there are
$a,b,c,d\in K$ such that $f|_{\langle u,v\rangle \times \langle u,v\rangle}$
is given by
\begin{equation} \eqlabel{Block}
    \begin{array}{c|cccc}
        f&1&u&v&uv\\
        \hline
        1&0&0&0&0\\
        u&0&a&a+b-c-d&-b+c+d\\
        v&0&d&b&b-d\\
        uv&0&a-d&-a+c+d&c
    \end{array}
\end{equation}
\end{proposition}
\begin{proof}
Let $u$, $v\in Q$ be such that $1\ne u\ne v\ne 1$. Then $\langle u,v\rangle =
\{1,u,v,w\}$, where $w=uv$. Our task is to define $f|_{\langle u,v\rangle
\times \langle u,v\rangle}$ so that it satisfies \peqref{MK1},
\peqref{MK2}. Since we must have $f(1,x)=f(x,1)=0$ by \peqref{cfs1}, it
remains to determine the nine entries $f(u,u)$, $f(u,v)$, $\dots$, $f(v,w)$,
$f(w,w)$. We think of the entries as variables.

Condition \peqref{MK1} yields $6$ equations in these variables, but only
three of these equations are distinct, namely
\begin{align*}
    f(u,v)&=f(v,v)-f(w,v),\\
    f(u,w)&=f(w,w)-f(v,w),\\
    f(v,u)&=f(u,u)-f(w,u).
\end{align*}
Similarly, condition \peqref{MK2} yields additional three equations
\begin{align*}
    f(u,v)&=f(u,u)-f(u,w),\\
    f(v,u)&=f(v,v)-f(v,w),\\
    f(w,u)&=f(w,w)-f(w,v).
\end{align*}
When we reorder the variables as $f(u,u)$, $f(v,v)$, $f(w,w)$, $f(u,v)$,
$f(u,w)$, $f(v,u)$, $f(v,w)$, $f(w,u)$, $f(w,v)$, the above equations
correspond to the system of linear equations
$A\overrightarrow{x}=\overrightarrow{0}$, where
\begin{displaymath}
    A=\left(\begin{array}{rrrrrrrrr}
        0&1&0&-1&0&0&0&0&-1\\
        0&0&1&0&-1&0&-1&0&0\\
        1&0&0&0&0&-1&0&-1&0\\
        1&0&0&-1&-1&0&0&0&0\\
        0&1&0&0&0&-1&-1&0&0\\
        0&0&1&0&0&0&0&-1&-1
    \end{array}
    \right).
\end{displaymath}
By keeping the 3rd, 5th and 6th row of $A$, it is now easy to see that by
subtracting rows from each other and by multiplying rows by $-1$, the system
$A$ is equivalent to
\begin{displaymath}
    \left(\begin{array}{rrrrrrrrr}
        1&0&0&0&0&-1&0&-1&0\\
        0&1&0&0&0&-1&-1&0&0\\
        0&0&1&0&0&0&0&-1&-1\\
        0&0&0&1&0&-1&-1&0&1\\
        0&0&0&0&1&0&1&-1&-1
    \end{array}\right).
\end{displaymath}
Therefore the variables $f(u,u)=a$, $f(v,v)=b$, $f(w,w)=c$, and $f(v,u)=d$
can be chosen arbitrarily, and the remaining variables are as claimed in
\peqref{Block}.
\end{proof}

Let us call a $4\times 4$ array with entries from an abelian group $K$ a
\emph{block} (\emph{over $K$}) if it is of the form \peqref{Block} for some
$a,b,c,d\in K$.

An arbitrary C-factor set $(\mathrm{id},f)$ with $Q$ Steiner can then be
obtained as follows: (i) assign $f(1,x)=f(x,1)=0$ for every $x\in Q$, (ii)
assign the diagonal entries $f(x,x)$ with $x\ne 1$ arbitrarily, (iii) if $f$
is not completed, assign arbitrarily any available $f(x,y)$, complete the
corresponding block $f|_{\langle x,y\rangle \times \langle x,y\rangle}$ and
repeat step (iii).

Note that this implies that for every block $B$ and for every Steiner loop
$Q$ there is a C-factor set $(\mathrm{id},f)$ such that $f|_{\langle
u,v\rangle \times \langle u,v\rangle}=B$ for some $u$, $v\in Q$.

If one block $B$ is fixed, and when for every $u$, $v\in Q$ with $1\ne u\ne
v\ne 1$ it is possible to premute $u$, $v$, $uv$ so that $f|_{\langle
u,v\rangle \times \langle u,v\rangle}=B$, we say that the C-factor set
$(\mathrm{id},f)$ is \emph{based on the block} $B$.

After a short reflection we see that:

\begin{lemma}\lemlabel{DiagB}
Let $B$ be a block with diagonal entries $a,b,c\in K$ such that
$a=b=c$. Then there is a C-factor set $(\mathrm{id},f)$ based on $B$.
\end{lemma}

For instance, the two blocks
\begin{displaymath}
     B_1(a) =
     \begin{array}{c|rrrr}
        f&1&u&v&uv\\
        \hline
        1&0&0&0&0\\
        u&0&a&a&0\\
        v&0&0&a&a\\
        uv&0&a&0&a
    \end{array},
    \quad\quad
    B_2(d) =
    \begin{array}{c|rrrr}
        f&1&u&v&uv\\
        \hline
        1&0&0&0&0\\
        u&0&0&-d&d\\
        v&0&d&0&-d\\
        uv&0&-d&d&0
    \end{array}
\end{displaymath}
give rise to C-factor sets (not necessarily uniquely determined). We will
take advantage of these blocks in the next section.

When the diagonal entries $a,b,c$ of $B$ are not all the same, no
C-factor set $(\mathrm{id},f)$ is based on $B$ provided $Q$ is sufficiently
large:

\begin{lemma}
Let $B$ be a block with diagonal entries $a,b,c$ such
that $|\{a,b,c\}|\ge 2$. If $|Q|>4$ then no C-factor set $(\mathrm{id},f)$
is based on $B$.
\end{lemma}

\begin{proof}
Assume that $f$ is based on $B$. Let $u$, $v\in Q$, $1\ne u\ne v\ne 1$. Then
$f|_{\langle u,v\rangle \times \langle u,v\rangle}=B$, and we can assume
without loss of generality that $f(u,u)=a$, $f(v,v)=b$. Furthermore, we can
assume that $c$ occurs once among $a$, $b$, $c$. Since $|Q|>4$, there are
$x$, $y \in Q\setminus\langle u, v\rangle$ such that $xy=u$. Since $f$ is
based on $B$, we have $f(x,x)=b$, $f(y,y)=c$, say. Now, $f|_{\langle
v,x\rangle \times \langle v,x\rangle}$ contains $c$ at least twice among its
diagonal entries, showing that $f$ is not based on B after all.
\end{proof}

\begin{example}
Let us illustrate the Lemma with the block
\begin{displaymath}
     B_3(a) =
     \begin{array}{c|rrrr}
        f&1&u&v&uv\\
        \hline
        1&0&0&0&0\\
        u&0&0&0&0\\
        v&0&-a&0&a\\
        uv&0&a&0&a
    \end{array}
\end{displaymath}
and with the smallest nonassociative Steiner loop of order $10$:
\begin{displaymath}
\begin{array}{r|rrrrrrrrrr}
        Q&0&1&2&3&4&5&6&7&8&9\\
        \hline
        0&0& 1& 2& 3& 4& 5& 6& 7& 8& 9\\
        1&1& 0& 3& 2& 5& 4& 7& 6& 9& 8\\
        2&2& 3& 0& 1& 6& 9& 4& 8& 7& 5\\
        3&3& 2& 1& 0& 8& 7& 9& 5& 4& 6\\
        4&4& 5& 6& 8& 0& 1& 2& 9& 3& 7\\
        5&5& 4& 9& 7& 1& 0& 8& 3& 6& 2\\
        6&6& 7& 4& 9& 2& 8& 0& 1& 5& 3\\
        7&7& 6& 8& 5& 9& 3& 1& 0& 2& 4\\
        8&8& 9& 7& 4& 3& 6& 5& 2& 0& 1\\
        9&9& 8& 5& 6& 7& 2& 3& 4& 1& 0
\end{array}
\end{displaymath}
(We used the block $B_3(a)$ in \cite[Proposition 3.4]{CLoopsI}, without
calling it ``block'', to construct the smallest nonassociative noncommutative
C-loop. Furthermore, the multiplication table of the $10$-element Steiner
loop is obtained from that of \cite{CLoopsI} by interchanging the elements
$6$ and $9$. The current multiplication table is more suitable here, as we
will see.)

Consider the subloop $\{0,1,2,3\}$. Define $f|_{\langle 1,2 \rangle \times
\langle 1,2 \rangle}$ so that it coincides with $B_3(a)$. It is possible to
extend $f$ so that $f|_{\langle 1,k \rangle \times \langle 1,k \rangle}$
coincides with $B_3(a)$ for every $k>0$. The partial map $f$ then looks as
follows:
\begin{displaymath}
\begin{array}{r|rrrrrrrrrr}
        Q&0&1&2&3&4&5&6&7&8&9\\
        \hline
        0&0& 0& 0& 0& 0& 0& 0& 0& 0& 0\\
        1&0& a& 0& a& 0& a& 0& a& 0& a\\
        2&0& a& 0& -a& & & & & & \\
        3&0& 0& 0& 0& & & & & & \\
        4&0& a& & & 0& -a& & & & \\
        5&0& 0& & & 0& 0& & & & \\
        6&0& a& & & & & 0& -a& & \\
        7&0& 0& & & & & 0& 0& & \\
        8&0& a& & & & & & & 0& -a\\
        9&0& 0& & & & & & & 0& 0
\end{array}
\end{displaymath}
However, $f|_{\langle 2,4 \rangle \times \langle 2,4 \rangle}$ then shows
that $f$ is not based on $B_3(a)$, no matter how it is actually defined.
\end{example}

\section{Nonflexible C-loops with prescribed center and factor by center}
\seclabel{NN}

\noindent In this section, we determine for which abelian groups $K$ and
Steiner loops $Q$ there exist nonflexible (hence noncommutative) C-loops $C$
with center $K$ and factor $C/K=Q$. Much of the discussion applies to nuclear
extensions as well.

Let us first say more about flexibility in extensions of C-loops:

\begin{lemma}
\lemlabel{NewMK}
In a C-loop $C$ with central squares, the following conditions are
equivalent.
\begin{enumerate}
\item[(i)] $C$ is flexible.

\item[(ii)] $(xy)^2 = (yx)^2$ for all $x,y\in C$.

\item[(iii)] $x\mapsto x^3$ is an antiautomorphism.
\end{enumerate}
\end{lemma}

\begin{proof}
First, $x (xy)^2 = x^2 y\cdot xy = x^2 (y\cdot xy)$ and so $(xy)^2 = x
(y\cdot xy)$. Next, $x (yx)^2 = (yx)^2 x = yx\cdot yx^2 = (yx\cdot y)x^2 =
x^2 (yx\cdot y)$ and so $(yx)^2 = x (yx\cdot y)$. This shows the equivalence
of (i) and (ii).

Next, $(xy)^2 = (xy)^3 \cdot y^{-1} x^{-1}$ and $(yx)^2 = yx\cdot (y^2
y^{-1}\cdot x^2 x^{-1}) = y^3 x^3\cdot y^{-1} x^{-1}$. This shows the
equivalence of (ii) and (iii).
\end{proof}

\begin{corollary}
\corolabel{NewMK}
A C-loop of exponent $4$ with central squares is flexible.
\end{corollary}

\begin{proof}
As noted, C-loops are inverse property loops. In any
inverse property loop of exponent $n$, $x\mapsto x^{n-1} = x^{-1}$ is an
antiautomorphism, and so the result follows from Lemma \lemref{NewMK}.
\end{proof}

\begin{corollary}
\corolabel{F2}
A C-loop with nucleus of order $2$ is flexible.
\end{corollary}

\begin{proof}
By \cite[Proposition 2.4]{CLoopsI}, the nucleus of a C-loop
contains every square. In any loop with nucleus of order $2$,
the nucleus coincides with the center. Now apply
Corollary \cororef{NewMK}.
\end{proof}

\begin{corollary}
\corolabel{Flex}
Let $Q$ be a Steiner loop, $K$ an abelian group, and $(\mathrm{id},f)$ a
C-factor set. Then $K\ltimes_{\mathrm{id}}^f Q$ is flexible if and only if $2
f(x,y) = 2 f(y,x)$ for all $x,y\in Q$.
\end{corollary}

\begin{proof}
For $a,b\in K$, $x,y\in Q$, we compute
\[
    ((a,x)\ast (b,y))^2 = (2a + 2b + 2f(x,y) + f(xy,xy),\,1).
\]
Reversing the roles of $(a,x)$ and $(b,y)$, we see that condition (ii) of
Lemma \lemref{NewMK} is satisfied if and only if $2f(x,y) = 2f(y,x)$ for all
$x,y\in Q$.
\end{proof}

We are now ready to characterize the parameters $K=Z(C)$, $Q=C/K$ of
nonflexible C-loops.

\begin{lemma}
\lemlabel{NoZ}
Let $C$ be a nonflexible C-loop with center $K$, and let
$Q = C/K$. Then $|Q|>2$ and $K$ is not an elementary abelian $2$-group.
\end{lemma}

\begin{proof}
If $|Q|=1$, then $C$ is associative, hence flexible. If $|Q|=2$, then $C$
has nucleus of index $2$, which is impossible by \cite[Lemma 2.9]{CLoopsI}.

Note that $C=K\ext{\mathrm{id}}{f}Q$ for some C-factor set $(\mathrm{id},f)$
by Lemma \lemref{Z}. If $K$ is an elementary abelian $2$-group, then
$2f(x,y)=0$ for every $x,y\in Q$. Then $C$ is flexible by Corollary \cororef{Flex}.
\end{proof}

\begin{lemma}
\lemlabel{Constr1}
Let $Q$ be a Steiner loop with $|Q|>2$, and let $K$ be an abelian
group which is not an elementary abelian $2$-group. Then there exists a
nonflexible C-loop $C$ of order $|Q|\cdot|K|$ such that the
nucleus $N = N(C)$ is isomorphic to $K$, and such that $C/N$ is
isomorphic to $Q$. Moreover, $N = Z(C)$.
\end{lemma}

\begin{proof}
Let $a\in K$ be an element of order different from $1$ or $2$. By Lemma
\lemref{DiagB}, there is a C-factor set $(\mathrm{id},f)$ with $f:Q\times
Q\to K$ based on the block $B=B_1(a)$. Let $C=K\ext{\mathrm{id}}{f}Q$ be the
corresponding C-loop.
Given $1\ne x\in Q$, we see from $B$ that there is $y\in Q$ such that
$f(x,y)=a$ and $f(y,x)=0$. Corollary \cororef{Flex} then shows that $C$ is not
flexible. In fact, a quick calculation yields $(b,x)(c,y)\cdot (b,x)\ne
(b,x)\cdot (c,y)(b,x)$, where we choose $y$ as above, and where $b,c\in K$ are
arbitrary. Since $N\le K$ by Theorem \thmref{ext}, $N=K$ follows.
By Lemma \lemref{Z}, $N=Z(C)$.
\end{proof}

In summary:

\begin{theorem}
\thmlabel{NN}
A nonflexible C-loop $C$ with center $K$ and Steiner factor
$C/K=Q$ exists if and only if $Q$ is a Steiner loop with $|Q|>2$
and $K$ is not an elementary abelian $2$-group. Moreover, it is
possible to demand $Z(C)=N(C)$ in such a case.
\end{theorem}

\subsection{Nonflexible C-loops with prescribed abelian
nucleus and factor by nucleus}
\subseclabel{NN}

\noindent We now turn to the following related question:

\begin{problem}
\problabel{1}
For which abelian groups $K$ and Steiner loops $Q$ does there
exist a nonflexible C-loop $C$ with nucleus $K$ and $C/K=Q$?
\end{problem}

Theorem \thmref{NN} gives a partial answer. Corollary \cororef{F2} implies
that $|K|>2$. Lemma 2.9 of \cite{CLoopsI} implies that $|Q|>2$. Thus the only
situation that remains to be considered is: $K$ an elementary abelian
$2$-group of order greater than $2$, and $Q$ a Steiner loop of order greater
than $2$.

\begin{lemma}
\lemlabel{NKlein}
Let $Q$ be the Klein group and $K$ an elementary abelian
$2$-group with $|K|>2$. Then there exists a nonflexible
C-loop $C$ such that $N(C)=K$ and $C/K=Q$.
\end{lemma}

\begin{proof}
Let $Q=\{1,u,v,w\}$. Assume first that $K=Q$. (We thus return to the
multiplicative notation for $K$ in this proof.) Define $\theta:Q\to\aut{K}$
by $\theta_1=\theta_u=\mathrm{id}_K$, $\theta_v=\theta_w=(v,w)$, where
$(v,w)$ is the transposition of $v$ and $w$. Define $f:Q\times Q\to K$ by
\[
    \begin{array}{c|cccc}
        f&1&u&v&w\\
        \hline
        1&1&1&1&1\\
        u&1&1&1&1\\
        v&1&v&1&w\\
        w&1&v&1&v
    \end{array}.
\]
A straightforward but tedious verification of condition \peqref{cfs2}
then shows that $(\theta,f)$ is a C-factor set. Moreover,
\begin{align*}
    &(u,v)(v,u)=(u\theta_vvf(v,u),vu)=(uwv,w)=(1,w)\\
    &\ne (w,w) = (vu,w) =(v\theta_uuf(u,v),uv) = (v,u)(u,v).
\end{align*}
and consequently
\begin{align*}
    &(u,v)(v,u)\cdot (u,v) = (1,w)(u,v) = (\theta_wuf(w,u),wv)=(uv,u)=(w,u)\\
    &\ne (1,u) = (uvw,u) = (u\theta_vwf(v,w),vw) = (u,v)(w,w) = (u,v)\cdot
    (v,u)(u,v).
\end{align*}
Thus $C_{16}=K\ext{\theta}{f}Q$ is a nonflexible C-loop with
the desired nucleus and nuclear factor.

To obtain a C-loop whose nucleus is a given elementary abelian $2$-group $K$,
$|K|=2^m>4$ and with nuclear factor $Q$, it suffices to take the direct
product of the loop $C_{16}$ with the elementary abelian $2$-group of order
$2^{m-2}$.
\end{proof}

Table \tableref{16} gives a multiplication
table of one of smallest nonflexible noncommutative C-loops with nucleus
isomorphic to the Klein group. The loop is constructed via the extension
described in the proof of Lemma \lemref{NKlein}.

\begin{table}[htb]
\caption{A smallest nonflexible noncommutative C-loop with nucleus that is an
elementary abelian $2$-group}
\[
\begin{array}{c||cccc|cccc|cccc|cccc}
   &0& 1& 2& 3& 4& 5& 6& 7& 8& 9& 10& 11& 12& 13& 14& 15 \\
   \hline\hline
   0&0& 1& 2& 3& 4& 5& 6& 7& 8& 9& 10& 11& 12& 13& 14& 15 \\
   1&1& 0& 3& 2& 5& 4& 7& 6& 9& 8& 11& 10& 13& 12& 15& 14 \\
   2&2& 3& 0& 1& 6& 7& 4& 5& 10& 11& 8& 9& 14& 15& 12& 13 \\
   3&3& 2& 1& 0& 7& 6& 5& 4& 11& 10& 9& 8& 15& 14& 13& 12 \\
   \hline
   4&4& 5& 7& 6& 0& 1& 3& 2& 12& 13& 15& 14& 8& 9& 11& 10 \\
   5&5& 4& 6& 7& 1& 0& 2& 3& 13& 12& 14& 15& 9& 8& 10& 11 \\
   6&6& 7& 5& 4& 2& 3& 1& 0& 14& 15& 13& 12& 10& 11& 9& 8 \\
   7&7& 6& 4& 5& 3& 2& 0& 1& 15& 14& 12& 13& 11& 10& 8& 9 \\
   \hline
   8&8& 9& 11& 10& 14& 15& 13& 12& 0& 1& 3& 2& 7& 6& 4& 5 \\
   9&9& 8& 10& 11& 15& 14& 12& 13& 1& 0& 2& 3& 6& 7& 5& 4 \\
   10&10& 11& 9& 8& 12& 13& 15& 14& 2& 3& 1& 0& 5& 4& 6& 7 \\
   11&11& 10& 8& 9& 13& 12& 14& 15& 3& 2& 0& 1& 4& 5& 7& 6 \\
   \hline
   12&12& 13& 14& 15& 10& 11& 8& 9& 4& 5& 6& 7& 2& 3& 0& 1 \\
   13&13& 12& 15& 14& 11& 10& 9& 8& 5& 4& 7& 6& 3& 2& 1& 0 \\
   14&14& 15& 12& 13& 8& 9& 10& 11& 6& 7& 4& 5& 0& 1& 2& 3 \\
   15&15& 14& 13& 12& 9& 8& 11& 10& 7& 6& 5& 4& 1& 0& 3& 2
\end{array}
\]
\tablelabel{16}
\end{table}

We have not been able to fully answer Problem \probref{1}. Note, however,
that it is not possible to prescribe an arbitrary quotient $Q$, since Lemma
\lemref{sim} and Theorem \thmref{NN} imply:

\begin{corollary}
Let $K$ be an elementary abelian $2$-group and $Q$ a simple Steiner loop.
Then there is no nonflexible C-loop $C$ with nucleus $K$ and
$C/K=Q$.
\end{corollary}

\subsection{Associators of given orders}
\subseclabel{ass}

\noindent Consider \cite[Proposition 3.4]{CLoopsI} restated in terms of
blocks:

\begin{proposition}
\proplabel{Constr}
Let $Q=\{1,u,v,uv\}$ be the Klein group. Let $K$ be an abelian group and
$a\in K$ an element of order at least $3$. Assume that $f:Q\times Q\to K$ is
based on the block $B_3(a)$. Then $K\ext{\mathrm{id}}{f}Q$ is a nonflexible
C-loop with nucleus and center $Z=N=\{(a,1);\;a\in K\}$.
\end{proposition}

The extension of Proposition \propref{Constr} allows us to construct a
nonflexible C-loop with an associator of given order $n>2$.

Let $n>2$, $\langle a\rangle=K=(\mathbb Z_n,+)$, $Q=\{1,u,v,w\}$ and $f$ be
as in Proposition \propref{Constr}. Let $x=(0,w)$, $y=(0,u)$, $z=(a,w)\in
C=K\ext{\mathrm{id}}{f}Q$. Then $((xy)z)^{-1} = ((0,w)(0,u)\cdot(a,w))^{-1} =
((a,v)(a,w))^{-1}=(3a,u)^{-1} = (-3a,u)$. Also, $x\cdot yz =
(0,w)\cdot(0,u)(a,w) = (0,w)(a,v)=(a,u)$. Therefore the associator of $x$,
$y$, $z$ is equal to $[x,y,z]=(xy\cdot z)^{-1}\cdot(x\cdot yz) =
(-3a,u)(a,u)=(-2a,1)$. It follows that the order of $[x,y,z]$ is $n$ when $n$
is odd, and $n/2$ when $n$ is even.

\begin{corollary}
\corolabel{assoc}
For every $n>2$ there is a nonflexible C-loop
$C$ containing an associator of order $n$. Moreover, it is possible to have
$|C|=4n$ if $n$ is odd, and $|C|=8n$ if $n$ is even.
\end{corollary}

Corollary \cororef{assoc} should be contrasted with the case of
extra loops. In an extra loop, every associator has order $2$
\cite[Theorem 5.7]{KKP}.

We have some evidence for the following:

\begin{conjecture}
Let $C$ be a C-loop with an associator of order $n$. Then $|C|\ge 4n$ if
$n$ is odd, and $|C|\ge 8n$ if $n$ is even.
\end{conjecture}

\section{The standard Cayley-Dickson process}
\seclabel{cd}

\noindent In view of Proposition \propref{Constr} and Corollary \cororef{F2},
it is now natural to ask if there is a flexible, nonextra, noncommutative
C-loop with nucleus (hence center) of order $2$. The answer is ``yes'',
as we will see in \S\secref{sed}, and the loop can be obtained in a classical
way---by the standard Cayley-Dickson process, which we deliberately define
in a broad way here.

The notation and terminology of this section are taken mostly from
\cite{SpVe}.

By an \emph{algebra} we mean here a vector space (over a field $F$) with
multiplication that distributes over addition, and with a neutral element $1$
with respect to multiplication. In particular, we do not assume that the
multiplication is associative.

\subsection{The process}
\subseclabel{process}

Let $D$ be a finite-dimensional algebra over $F$, let $N_D:D\to F$ be a
quadratic form, and set $\lambda=-1\in F$. Let $\langle\ ,\ \rangle_D:D\oplus
D\to F$ be the associated bilinear form defined by $\langle x,y\rangle_D =
N_D(x+y) - N_D(x) - N_D(y)$, and let $\ov{x}=\langle x,1\rangle_D 1-x$ be the
\emph{conjugate} of $x\in D$.

We say that $(C,N_C)$ is obtained from $(D,N_D)$ by \emph{standard doubling}
if
\begin{enumerate}
\item[(i)] $C=D\oplus D$ as a vector space,

\item[(ii)] the multiplication in $C$ is given by
\begin{equation} \eqlabel{Doubling}
    (x,y)(u,v)=(xu+\lambda\ov{v}y, vx+y\ov{u}),
\end{equation}

\item[(iii)] $N_C((x,y)) = N_D(x)-\lambda N_D(y)$.
\end{enumerate}

Since $\langle\ ,\ \rangle_D$ is bilinear, it follows that conjugation is an
additive homomorphism, and it is easy to check that $C$ is an algebra and
that $N_C$ is a quadratic form. It is therefore possible to iterate the
standard doubling. It is customary to call these iterations the
\emph{standard Cayley-Dickson process}.

\subsection{Structure constants}
\subseclabel{constants}

The presence of the distributive laws makes it possible and convenient to
describe the multiplication in $D$ by \emph{structure constants}. In particular,
when $\{d_1,\dots,d_n\}$ is a basis for $D$, then there are $n^3$ constants
$\gamma_{ij}^k\in F$ such that $d_id_j=\sum_{k=1}^n \gamma_{ij}^kd_k$.

By selecting the basis of $D$ carefully and systematically in the standard
doubling, all but one structure constants can be eliminated for given $i,j$,
as we are going to show. (The usual multiplication formulae for
quaternions and octonions are based on this observation.)

\begin{lemma}
\lemlabel{Conj}
Assume that $(C,N_C)$ is obtained from $(D,N_D)$ by one application of
doubling. Then $\ov{(u,v)}=(\ov{u},-v)$ for every $u$, $v\in D$.
\end{lemma}

\begin{proof}
We have $\ov{(u,v)} = \langle(u,v),(1,0)\rangle_C(1,0) - (u,v) = [
N_C((u,v)+(1,0)) - N_C(u,v) - N_C(1,0)](1,0) - (u,v) = [N_D(u+1) + N_D(v) -
N_D(u) - N_D(v) - N_D(1) - N_D(0)](1,0) - (u,v) = \langle u,1\rangle_D(1,0) -
(u,v) = (\ov{u}+u)(1,0) - (u,v) = (\ov{u},-v)$.
\end{proof}

\begin{lemma}
\lemlabel{SimpleBasis}
Let $D$ be an algebra over $F$, $N_D:D\to F$ a quadratic form, and let
$(C,N_C)$ be obtained from $(D,N_D)$ by the standard doubling. Suppose that
$\{d_1,\dots,d_n\}$ is a basis of $D$ such that:
\begin{enumerate}
\item[(i)] $d_1=1$,

\item[(ii)] $\ov{d_i}=-d_i$ for $1<i\le n$,

\item[(iii)] $d_id_j\in\pm\{d_1,\dots,d_n\}$ for $1\le i$, $j\le n$,

\item[(iv)] $d_id_i=-1$ for $1<i\le n$,

\item[(v)] $d_id_j=-d_jd_i$ for $1<i<j\le n$,

\item[(vi)] $\ov{d_id_j} = \ov{d_j}\cdot\ov{d_i}$, $\ov{\ov{d_i}}=d_i$, for
$1\le i$, $j\le n$,

\item[(vii)] $d_i(d_id_j)= (d_id_i)d_j$, $d_i(d_jd_j)=(d_id_j)d_j$,
$(d_id_j)d_i = d_i(d_jd_i)$ for $1\le i$, $j\le n$.
\end{enumerate}
For $1\le i\le n$, let $c_i=(d_i,0)$, $c_{i+n}=(0,d_i)$. Then
$\{c_1,\dots,c_{2n}\}$ is a basis of $C$ satisfying the properties analogous
to $(\mathrm{i})$--$(\mathrm{vii})$, with $2n$ instead of $n$.
\end{lemma}

\begin{proof}
This follows from straightforward calculation, but we show most of the
proof for completeness.

Part (i) is trivial. By Lemma \lemref{Conj}, $\ov{(d_i,0)}=
(\ov{d_i},0)=(-di,0)$, $\ov{(0,d_i)} = (0,-d_i)$ for every $1<i\le n$, and
(ii) follows.

While multiplying $c_i$ by $c_j$, we notice that only one of the four
summands in \peqref{Doubling} is nonzero, and this nonzero summand has
coefficient $\pm 1$. Namely, for $1\le i$, $j\le n$, we have
\[
    c_ic_j=(d_id_j,0),\quad c_ic_{j+n}=(0,d_jd_i),\quad
    c_{i+n}c_j=(0,d_i\ov{d_j}),\quad c_{i+n}c_{j+n}=(-\ov{d_j}d_i,0).
\]
This shows (iii).

Part (iv) is easy, and we proceed to prove (v). With $1<i<j\le n$ we get
$c_ic_j = (d_id_j,0) = (-d_jd_i,0) = -c_jc_i$. With $1<i\le n$, $1\le j\le n$
we get $c_ic_{j+n} = (0,d_jd_i) = -(0,d_j\ov{d_i})=-c_{j+n}c_i$. When $1\le
i<j\le n$, we have $c_{i+n}c_{j+n} = (-\ov{d_j}d_i,0) = -(-\ov{d_i}d_j,0) =
-c_{j+n}c_{i+n}$.

Part (vi) is similar, and can be simplified by noting that it suffices to
prove (vi) for $i<j$.

We now prove $c_i(c_ic_j)=(c_ic_i)c_j$ for $1\le i$, $j\le 2n$, and leave the
other two parts of (vii) to the reader. There is nothing to show when $i=1$
or $j=1$. When $i>1$, the equation becomes $c_i(c_ic_j)=-c_j$. With $1<i$,
$j\le n$, we then get $c_i(c_ic_j) = (d_i(d_id_j),0) = (-d_j,0) = -c_j$,
$c_i(c_ic_{j+n}) = (d_i,0)(0,d_jd_i) = (0,(d_jd_i)d_i) = (0, -d_j) =
-c_{j+n}$, $c_{i+n}(c_{i+n}c_{j+n}) = (0,d_i)(-\ov{d_j}d_i,0) =
(0,d_i\ov{(-\ov{d_j}d_i)}) = -(0,d_i(\ov{d_i}d_j)) = -(0,d_j) = -c_{j+n}$,
where in the last case we allow $i=1$.
\end{proof}

\subsection{Remarks on the Cayley-Dickson process}
\subseclabel{remarks}

The word \emph{standard} in standard Cayley-Dickson process refers to the
fact that $\lambda=-1$. The process makes sense for any value of $\lambda\in
F^*$.

It is usually assumed that $N_D$ is a nondegenerate quadratic form satisfying
$N_D(uv)=N_D(u)N_D(v)$. Then $(D,N_D)$ is called a \emph{composition
algebra}.

When the underlying field $F$ is the real numbers $\mathbb R$, we speak of
\emph{real algebras}. Nevertheless, the characteristic of $F$ can be
arbitrary.

The classical Hurwitz theorem has been extended to all characteristics and
composition algebras, i.e., composition algebras exist only in dimension $1$,
$2$, $4$ and $8$. See \cite[Chapter 1]{SpVe} for details and \cite[Theorem
1.6.2]{SpVe} in particular.

Smith \cite{Sm} constructed a real 16-dimensional composition semialgebra,
which does not satisfy one of the distributive laws but otherwise has all the
properties of a composition algebra. Kivunge and Smith \cite{kvs} study
subloops (which are not C-loops unless associative) of the associated
left loop.

\section{C-loops arising from the standard Cayley-Dickson process}
\seclabel{sed}

\noindent Let $A_n$ be the algebra of dimension $2^n$ constructed from
$A_0=\mathbb R$ with $N_\mathbb R:x\mapsto x^2$ by the standard
Cayley-Dickson process. Assume that the basis $\{a_1,\dots,a_{2^n}\}$ of
$A_n$ is obtained systematically as in Lemma \lemref{SimpleBasis}. Then
$L_n=\pm\{a_1,\dots,a_{2^n}\}$ is closed under multiplication by Lemma
\lemref{SimpleBasis}(iii), and it possesses a neutral element $1=a_1$. Note
that $-a_1$ commutes and associates with all elements of $L_n$. Also note
that the inverse of $a_i$ is either $a_i$ or $-a_i$, and that $L_n$ satisfies
the alternative and flexible laws, by Lemma \lemref{SimpleBasis}(vii). Hence
$L_n$ is a flexible, alternative, inverse property loop.

\begin{proposition}
\proplabel{FC}
For every $n\ge 1$, $L_n$ is a flexible C-loop.
\end{proposition}

\begin{proof}
Let $x$, $y$, $z\in L_n$. Then $x(y(yz))=x(y^2z)$ by alternativity. By Lemma
\lemref{SimpleBasis}, $yy=\pm 1$, and, as we have observed, $\pm 1 \in
Z(L_n)$. Therefore $x(y^2 z) = (xz)y^2$. Similarly, $((xy)y)z= (xz)y^2$.
\end{proof}

\begin{proposition}
\proplabel{COMP}
For $n\ge 1$, let $L_n$ be the C-loop of signed basic elements of $A_n$.
Then:
\begin{enumerate}
\item[(i)] $N(L_2)=Z(L_2)=L_2$,

\item[(ii)] $N(L_3)=L_3$, $Z(L_3) = \{\pm 1\}$,

\item[(iii)] $N(L_n)=Z(L_n)=\{\pm 1\}$ for every $n>0$, $n\not\in\{2,3\}$.
\end{enumerate}
\end{proposition}

\begin{proof}
Using the standard notation for complex numbers and quaternions, we have
$L_1=\{\pm 1\}$, $L_2=\{\pm 1,\pm i\}\cong \mathbb Z_4$, $L_3=\{\pm 1, \pm i,
\pm j, \pm k\}\cong Q$, where $Q$ is the quaternion group. It remains to show
(iii).

We will demonstrate later (outside of this proof) that (iii) holds for $L_4$.
Assume that (iii) holds for every $L_m$ with $4\le m\le n$.

If $x=(u,0)$ for $u\in L_n\setminus\{\pm 1\}$ then $x\not\in N(L_{n})$ by
induction assumption, and since $L_{n}\le L_{n+1}$, we have $x\not\in
N(L_{n+1})$.

Assume that $x=(0,u)$ for $u\in L_n$. For $v,w\in L_n$, we have
\begin{align} \eqlabel{Aux1}
    &(0,1)(v,0)\cdot (0,w) =
    (0,\ov{v})(0,w)=(-\ov{w}\cdot\ov{v},0),\\
     \eqlabel{Aux2}
    &(0,1)\cdot(v,0)(0,w)=(0,1)(0,wv)=(-\ov{wv},0).
\end{align}
When $u=\pm 1$, pick $v,w\in L_n$ that do not commute (which is possible
by Lemma \lemref{SimpleBasis}), and conclude from $\peqref{Aux1})$,
$\peqref{Aux2})$ that $(0,u)\not\in N(L_{n+1})$. When $u\ne\pm 1$, let
$u=w$, pick $v\in L_n$ that does not commute with $u$, and conclude as in the
previous case that $(0,u)\not\in N(L_{n+1})$.
\end{proof}

\subsection{The standard sedenion loop}
\subseclabel{sedloop}

The C-loop $S=L_4=\pm\{a_1,\dots,a_{16}\}$ is called the \emph{standard
sedenion loop}. The ``structure
constants'' $a_ia_j=\gamma_{ij}$ of the $16$-dimensional algebra $A_4$ are as
in Table \tableref{SC}, where $i$ stands for $a_i$, and $-i$ for $-a_i$.

\begin{table}[htb]
\caption{Structure constants of the standard real sedenions $S$}
\begin{tiny}
\[
\begin{array}{r|rrrrrrrrrrrrrrrr}
\gamma_{ij}&1&2&3&4&5&6&7&8&9&19&11&12&13&14&15&16\\
\hline
1&1& 2& 3& 4& 5& 6& 7& 8& 9& 10& 11& 12& 13& 14& 15& 16 \\
2&2& -1& 4& -3& 6& -5& -8& 7& 10& -9& -12& 11& -14& 13& 16& -15 \\
3&3& -4& -1& 2& 7& 8& -5& -6& 11& 12& -9& -10& -15& -16& 13& 14 \\
4&4& 3& -2& -1& 8& -7& 6& -5& 12& -11& 10& -9& -16& 15& -14& 13 \\
5&5& -6& -7& -8& -1& 2& 3& 4& 13& 14& 15& 16& -9& -10& -11& -12 \\
6&6& 5& -8& 7& -2& -1& -4& 3& 14& -13& 16& -15& 10& -9& 12& -11 \\
7&7& 8& 5& -6& -3& 4& -1& -2& 15& -16& -13& 14& 11& -12& -9& 10 \\
8&8& -7& 6& 5& -4& -3& 2& -1& 16& 15& -14& -13& 12& 11& -10& -9 \\
9&9& -10& -11& -12& -13& -14& -15& -16& -1& 2& 3& 4& 5& 6& 7& 8 \\
10&10& 9& -12& 11& -14& 13& 16& -15& -2& -1& -4& 3& -6& 5& 8& -7 \\
11&11& 12& 9& -10& -15& -16& 13& 14& -3& 4& -1& -2& -7& -8& 5& 6 \\
12&12& -11& 10& 9& -16& 15& -14& 13& -4& -3& 2& -1& -8& 7& -6& 5 \\
13&13& 14& 15& 16& 9& -10& -11& -12& -5& 6& 7& 8& -1& -2& -3& -4 \\
14&14& -13& 16& -15& 10& 9& 12& -11& -6& -5& 8& -7& 2& -1& 4& -3 \\
15&15& -16& -13& 14& 11& -12& 9& 10& -7& -8& -5& 6& 3& -4& -1& 2 \\
16&16& 15& -14& -13& 12& 11& -10& 9& -8& 7& -6& -5& 4& 3& -2& -1
\end{array}
\]
\end{tiny}
\tablelabel{SC}
\end{table}

The loop $S$ is then easily obtained. By Proposition \propref{FC},
$S$ is a flexible C-loop. One can check by hand (or by computer) that
$N(S)=Z(S)=\{\pm1\}$, thus completing the proof of Proposition
\propref{COMP}.

Moreover, $S$ is not extra. Therefore any of the loops $L_n$, $n>3$, is a
flexible, nonextra, noncommutative C-loop with nucleus of order $2$. These
are the loops we set out to find at the beginning of \S\secref{cd}.

More detailed information about the standard sedenion loop,
including its subloop structure, can be found in \cite{Cawagas}.
The loop $S$ is contained in the library of loops
of the GAP \cite{GAP4} package LOOPS \cite{LOOPS}. The above-mentioned
properties of $S$ can be verified easily with LOOPS.

%%%%%%%%%%%%%%%%%%%%%%%%%%%
% BIBLIOGRAPHY
%%%%%%%%%%%%%%%%%%%%%%%%%%%

\bibliographystyle{plain}

\end{document}